\documentclass{amsart}
\usepackage{psfig}

\newtheorem{theorem}{Theorem}[section]
\newtheorem{lemma}[theorem]{Lemma}

\theoremstyle{definition}

\theoremstyle{remark}

\numberwithin{equation}{section}

\newcommand{\abs}[1]{\lvert#1\rvert}

\DeclareSymbolFont{AMSb}{U}{msb}{m}{n}
\DeclareMathSymbol{\Z}{\mathalpha}{AMSb}{"5A}

\begin{document}
\newcommand{\beqs}{\begin{equation*}}
\newcommand{\eeqs}{\end{equation*}}
\newcommand{\beq}{\begin{equation}}
\newcommand{\eeq}{\end{equation}}
\newcommand\nutwid{\overset {\text{\lower 3pt\hbox{$\sim$}}}\nu}
\newcommand\Mtwid{\overset {\text{\lower 3pt\hbox{$\sim$}}}M}
\newcommand\bijone{\overset {1}\longrightarrow}                   
\newcommand\bijtwo{\overset {2}\longrightarrow}                   
\newcommand\pihat{\widehat{\pi}}
\newcommand\mymod[1]{(\mbox{mod}\ {#1})}

\title[The Andrews-Stanley refinement of Ramanujan's congruence]{On the Andrews-Stanley Refinement of \\
Ramanujan's Partition Congruence Modulo $5$}

\author{Alexander Berkovich}
\address{Department of Mathematics, University of Florida, Gainesville,
Florida 32611-8105}
\email{alexb@math.ufl.edu}          

\author{Frank G. Garvan}
\address{Department of Mathematics, University of Florida, Gainesville,
Florida 32611-8105}
\email{frank@math.ufl.edu}          


\date{December 26, 2003} 



\begin{abstract}
In a recent study of sign-balanced, labelled posets Stanley \cite{Stan}, 
introduced a new integral partition statistic
\beqs
\mathrm{srank}(\pi) = {\mathcal O}(\pi) - {\mathcal O}(\pi'),
\eeqs
where ${\mathcal O}(\pi)$ denotes the number of odd parts of the partition $\pi$
and $\pi'$ is the conjugate of $\pi$. In \cite{Andrews1} Andrews proved
the following refinement of Ramanujan's partition congruence mod $5$:
\begin{align*}
p_0(5n+4) &\equiv p_2(5n+4) \equiv 0 \pmod{5},\\
p(n) &= p_0(n) + p_2(n),
\end{align*}
where $p_i(n)$ ($i=0,2$) denotes the number of partitions of $n$ with
$\mathrm{srank}\equiv i\pmod{4}$ and $p(n)$ is the number of unrestricted partitions
of $n$. Andrews asked for a partition statistic that would divide the
partitions enumerated by $p_i(5n+4)$ ($i=0,2$) into five equinumerous classes.

In this paper we discuss two such statistics. The first one, while new, is 
intimately related to the Andrews-Garvan \cite{AG} crank. The second one is in terms of
the $5$-core crank, introduced by Garvan, Kim and Stanton \cite{GKS}.
Finally, we discuss some new formulas for partitions that are $5$-cores.
\end{abstract}

\maketitle

\section{Introduction} \label{sec:intro}
Let $p(n)$ be the number of unrestricted partitions of $n$. 
Ramanujan discovered and later proved that
\begin{align}
p(5n+4) &\equiv 0 \pmod{5},\label{eq:ram5} \\ 
p(7n+5) &\equiv 0 \pmod{7},\label{eq:ram7} \\ 
p(11n+6) &\equiv 0 \pmod{11}.\label{eq:ram11} 
\end{align}
Dyson \cite{Dyson} was the first to consider combinatorial explanations
of these congruences. He defined the rank of a partition as 
{\it the largest part minus the number of parts} and made the
empirical observations that
\begin{align}
N(k,5,5n+4) &= \frac{p(5n+4)}{5},\quad 0 \le k \le 4, \label{eq:dys5} \\ 
N(k,7,7n+5) &= \frac{p(7n+5)}{7},\quad 0 \le k \le 6, \label{eq:dys7}  
\end{align}
where $N(k,m,n)$ denotes the number of partitions of $n$ with
rank congruent to $k$ modulo $m$.
Equation (\ref{eq:dys5}) means that the residue of
the rank mod ${5}$ divides the partitions of $5n+4$ 
into five equal classes. Similarly, 
(\ref{eq:dys7}) implies that the residue of
the rank mod ${7}$ divides the partitions of $7n+5$
into seven equal classes. Dyson's rank failed to explain (\ref{eq:ram11}),
and so Dyson conjectured the existence of a hypothetical statistic, called the
crank, that would explain the Ramanujan congruence mod ${11}$.
Identities (\ref{eq:dys5})-(\ref{eq:dys7}) were later proved by
Atkin and Swinnerton-Dyer \cite{ASD}. Andrews and Garvan \cite{AG} found a
crank for all three Ramanujan congruences (\ref{eq:ram5})-(\ref{eq:ram11}).
Their crank is defined as follows
\begin{equation}
\mbox{crank}(\pi) =
\begin{cases}
  \ell(\pi), &\mbox{if $\mu(\pi)=0$}, \\
  \nutwid(\pi) - \mu(\pi), &\mbox{if $\mu(\pi)>0$},
\end{cases}
\label{eq:crankdef} 
\end{equation}
where $\ell(\pi)$ denotes the largest part of $\pi$,
$\mu(\pi)$ denotes the number of ones in $\pi$ and $\nutwid(\pi)$
denotes the number of parts of $\pi$ larger than $\mu(\pi)$.

Later, Garvan, Kim and Stanton \cite{GKS} found different cranks, which also
explained all three congruences (\ref{eq:ram5})-(\ref{eq:ram11}).
Their approach made essential use of $t$-cores of partitions 
and led to explicit bijections between various equinumerous classes.
In particular, they provided what amounts to a combinatorial proof
of the formula
\beq
\sum_{n\ge0} p(5n+4) q^n
= 5 \prod_{m\ge1} \frac{ (1-q^{5m})^5 }{ (1-q^m)^6 },
\label{eq:rambest} 
\eeq
considered by Hardy to be an example of Ramanujan's best work.

The main results of \cite{AG} can be summarized as
\begin{align}
M(k,5,5n+4) &= \frac{p(5n+4)}{5},\quad 0 \le k \le 4, \label{eq:ag5}\\ 
M(k,7,7n+5) &= \frac{p(7n+5)}{7},\quad 0 \le k \le 6, \label{eq:ag7}\\ 
M(k,11,11n+6) &= \frac{p(11n+6)}{11},\quad 0 \le k \le 10, \label{eq:ag11} 
\end{align}
and
\begin{align}
&1 + (x + x^{-1} - 1)q + \sum_{n>1} \sum_{m} \Mtwid(m,n) x^m q^n \nonumber\\
&\qquad\qquad= 
\prod_{n\ge1} \frac{ (1 - q^n) }{ (1 - x q^n) (1 - x^{-1} q^n) },
\label{eq:crankgf} 
\end{align}
where $\Mtwid(m,n)$ denotes the number of partitions of $n$ with crank $m$
and $M(k,m,n)$ denotes the number of partitions of $n$ with
crank congruent to $k$ modulo $m$.

In \cite{G2} Garvan found a refinement of (\ref{eq:ram5})
\beq
M(k,2,5n+4) \equiv 0 \pmod{5}, \quad k=0,1 \label{eq:gref5} 
\eeq
together with the combinatorial interpretation
\beq
M(2k+\alpha,10,5n+4) = \frac{M(\alpha,2,5n+4)}{5},\quad 0 \le k \le 4, 
\label{eq:gref5a}\\ 
\eeq
with $\alpha=0,1$.

Recently, a very different refinement of (\ref{eq:ram5}) was given by Andrews
\cite{Andrews1}. Building on the work of Stanley \cite{Stan}, Andrews examined
partitions $\pi$ classified according to ${\mathcal O}(\pi)$ and
${\mathcal O}(\pi')$, where 
where ${\mathcal O}(\pi)$ denotes the number of odd parts of the partition $\pi$
and $\pi'$ is the conjugate of $\pi$. He used recursive relations to show that
\beq
G(z,y,q) := \sum_{n,r,s\ge0} S(n,r,s) q^n z^r y^s
=
\frac{ (-zyq;q^2)_\infty }{ (q^4;q^4)_\infty (z^2q^2;q^4)_\infty  
(y^2q^2;q^4)_\infty},
\label{eq:rsgf} 
\eeq
where $S(n,r,s)$ denotes the number of partitions $\pi$ of $n$ with
${\mathcal O}(\pi)=r$,
${\mathcal O}(\pi')=s$, and
\begin{align}
(a;q)_\infty &= \lim_{n\to\infty} (a;q)_n, \label{eq:aqdef} \\ 
(a;q)_n = (a)_n &=
\begin{cases}
1, &\mbox{if $n=0$},\\
\prod_{j=0}^{n-1}(1-aq^j), &\mbox{if $n>0$.}
\end{cases}
\label{eq:aqndef} 
\end{align}
A direct combinatorial proof of (\ref{eq:rsgf}) was later given
by A. Sills \cite{Sills}, A.~J.~Yee \cite{Y} and C. Boulet \cite{Boulet}. 
Actually, C. Boulet proved a stronger version of (\ref{eq:rsgf})
with one extra parameter.
We define the Stanley rank of a partition
$\pi$ as
\beq
\mathrm{srank}(\pi) = {\mathcal O}(\pi) - {\mathcal O}(\pi').
\label{eq:srankdef} 
\eeq
It is easy to see that
\beq
\mathrm{srank}(\pi) \equiv 0 \pmod{2},
\label{eq:srankcong} 
\eeq
so that
\beq
p(n) = p_0(n) + p_2(n),
\label{eq:p02} 
\eeq
where $p_i(n)$ ($i=0,2$) denotes the number of partitions of $n$ with
$\mathrm{srank}\equiv i\pmod{4}$.
We note that (\ref{eq:rsgf}) with $z=y^{-1}=\sqrt{-1}$ immediately implies
the Stanley formula \cite[p.8]{Stan}
\beq
\sum_{n\ge0} (p_0(n) - p_2(n)) q^n
= \frac{ (-q;q^2)_\infty }{ (q^4;q^4)_\infty (-q^2;q^4)_\infty^2}.
\label{eq:p02prod} 
\eeq
Using (\ref{eq:ram5}), (\ref{eq:p02}) and (\ref{eq:p02prod}), Andrews
proved the following refinement of (\ref{eq:ram5})
\beq
p_0(5n+4) \equiv p_2(5n+4) \equiv 0 \pmod{5}.
\label{eq:andrefine} 
\eeq
His proof of (\ref{eq:andrefine}) was analytic and so at the end of \cite{Andrews1}
he posed the problem of finding a partition statistic that would give a 
combinatorial interpretation of (\ref{eq:andrefine}).
The object of this paper is to provide a solution to the Andrews problem.
It turns out that there are two distinct integral partition statistics,
whose residue mod $5$ split the partitions enumerated by
$p_i(5n+4)$ (with $i=0,2$) into five equal classes. The first statistic,
which we call the stcrank, is new. However, it is intimately related
to the Andrews-Garvan crank (\ref{eq:crankdef}).
Unexpectedly, the second statistic is the ``$5$-core crank'', introduced
by Garvan, Kim and Stanton \cite{GKS}.
This second statistic not only provides the desired combinatorial
interpretation, but it also provides a direct combinatorial proof of
(\ref{eq:andrefine}).

The rest of this paper is organized as follows. In Section 2 we define the stcrank
and show that is indeed, a statistic asked for in \cite{Andrews1}. In Section 3
we briefly review the development in \cite{GKS}. In Section 4 we state a
number of new formulas for partitions that are $5$-cores. We sketch the
``$5$-core crank'' proof of (\ref{eq:andrefine}) and we conclude with
some open problems.

\section{The stcrank}
We begin with some preliminaries about partitions and their conjugates.
A partition $\pi$ is a nonincreasing sequence
\beq
\pi = (\lambda_1, \lambda_2, \lambda_3, \dots)
\label{eq:pidef} 
\eeq
of nonnegative integers (parts)
\beq
\lambda_1 \ge \lambda_2 \ge \lambda_3 \ge \cdots.
\label{eq:lams} 
\eeq
The weight of $\pi$, denoted by $\abs{\pi}$ is the sum of parts
\beq
\abs{\pi} = \lambda_1 + \lambda_2 + \lambda_3 + \cdots.
\label{eq:piweight} 
\eeq
If $\abs{\pi}=n$, then we say that $\pi$ is a partition of $n$. Often it is 
convenient to use another notation for $\pi$
\beq
\pi = (1^{f_1}, 2^{f_2}, 3^{f_3}, \dots),
\label{eq:pidef2}
\eeq
which indicates the number of times each integer occurs as a part.
The number $f_i=f_i(\pi)$ is called the frequency of $i$ in $\pi$.
The conjugate of $\pi$ is the partition 
$\pi'=(\lambda_1', \lambda_2', \lambda_3', \dots)$ with
\begin{align}
\lambda_1' &= f_1 + f_2 + f_3 + f_4 + \cdots \nonumber\\
\lambda_2' &= f_2 + f_3 + f_4 + \cdots \label{eq:piconj}\\
\lambda_3' &= f_3 + f_4 + \cdots \nonumber\\
	  &\vdots \nonumber
\end{align}
Next, we discuss two bijections. The first one relates $\pi$ and bipartitions
$(\pi_1,\pi_2)$, where $\pi_2$ is a partition with no repeated even parts.

\noindent
{\bf Bijection 1}  
\beqs
\pi \bijone (\pi_1,\pi_2),
\eeqs
where
\begin{align*}
\pi &= (1^{f_1}, 2^{f_2}, 3^{f_3}, \dots),\\
\pi_1 &= (1^{\lfloor{f_2/2}\rfloor}, 2^{\lfloor{f_4/2}\rfloor}, 
          3^{\lfloor{f_6/2}\rfloor}, \dots),\\
\pi_2 &= (1^{f_1}, 2^{\left\{f_2\right\}}, 3^{f_3}, 4^{\left\{f_2\right\}},
\dots),
\end{align*}
$\lfloor{x}\rfloor$ is the largest integer $\le x$, and
\beqs
\left\{x\right\} = x - 2 \lfloor{x/2}\rfloor.
\eeqs

\noindent
Indeed, remove from $\pi$ the maximum even number of even parts.
The resulting partition is $\pi_2$, The removed even parts can be organized
into a new partition
$(2^{2\lfloor{f_2/2}\rfloor}, 4^{2\lfloor{f_4/2}\rfloor},
6^{2\lfloor{f_6/2}\rfloor}, \dots),$
which can easily be mapped onto $\pi_1$.
Clearly, we have
\begin{align}
\abs{\pi} &= 4\abs{\pi_1} + \abs{\pi_2}, \label{eq:bij1prop1}\\ 
\mathrm{srank}(\pi) &= \mathrm{srank}(\pi_2), \label{eq:bij1prop2} 
\end{align}
so that
\begin{align}
\sum_{\pi} q^{\abs{\pi}} y^{\mathrm{srank}(\pi)}
&= \sum_{\pi_1} q^{4\abs{\pi_1}}
\sum_{\pi_2} q^{\abs{\pi_2}} y^{\mathrm{srank}(\pi_2)} \nonumber\\
&=\frac{1}{(q^4;q^4)_\infty} 
\sum_{\pi_2} q^{\abs{\pi_2}} y^{\mathrm{srank}(\pi_2)}.
\label{eq:bij1gf} 
\end{align}
Comparing (\ref{eq:bij1gf}) and  (\ref{eq:rsgf}) with $zy=1$, we see that
\beq
\sum_{\pi_2} q^{\abs{\pi_2}} y^{\mathrm{srank}(\pi_2)}
= \frac{ (-q;q^2)_\infty }{ (y^2q^2;q^4)_\infty (q^2/y^2;q^4)_\infty },
\label{eq:srankprodid} 
\eeq
where the sum is over all partitions with no repeated even parts.

To describe our second bijection we require a few definitions.
We say that $\pi_A$ is a partition of type A iff 
$\pi_A \bijone ((1),\pi_2)$. We say that 
$\pi_B = (\lambda_1, \lambda_2, \lambda_3, \dots)$ is
a partition of type B iff either $\abs{\pi_B}\ne4$,
$\lambda_1-\lambda_2\ge2$,
$\lambda_1'-\lambda_2'\ge2$,
$\lambda_1 -2$ and $\lambda_2$ are not identical even integers and
$\pi_B$ has no repeated even parts, or $\pi_B=(3,1)$.
Obviously, $\pi_B \bijone ((0),\pi_B)$.
Our second bijection relates partitions of type A and B.

\noindent
{\bf Bijection 2}  
\beqs
\pi_A \bijtwo \pi_B,
\eeqs
where
\begin{align*}
\pi_A &= (1^{f_1}, 2^{f_2}, 3^{f_3}, \dots, m^{f_m}),\\
\pi_B &=
\begin{cases}
(1^{f_1+2}, 2^{f_2-2}, 3^{f_3}, 4^{f_4}, \dots, (m-1)^{f_{m-1}}, m^{f_m-1},
(m+2)^1), &\mbox{if $m>2$}, \\
(1^{f_1+2},4^1), &\mbox{if $m=2$, $f_2=3$}, \\
(1^{f_1+1},3^1), &\mbox{if $m=2$, $f_2=2$},
\end{cases}
\end{align*}
$m\ge2$, $f_2=2$, $3$, and $f_{2i}=0$, $1$ for $i>1$.

\begin{figure}
\centerline{\psfig{figure=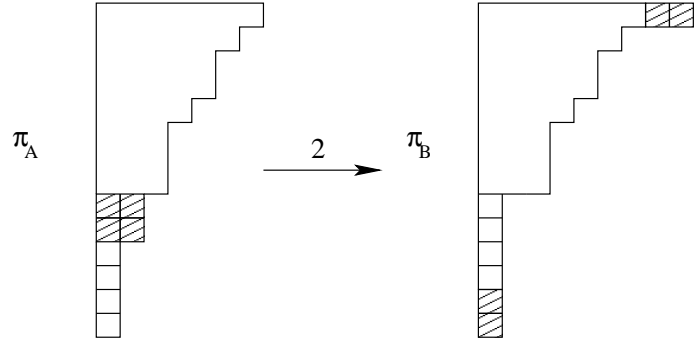}}
\caption{Graphical illustration of Bijection 2}
\label{fig1}
\end{figure}

Clearly, we have
\begin{align}
\abs{\pi_A} &= \abs{\pi_B}, \label{eq:bij2prop1} \\ 
\mathrm{srank}(\pi_A) &= \mathrm{srank}(\pi_B). \label{eq:bij2prop2} 
\end{align}
Next, we define a new partition statistic
\beq
\mathrm{stcrank}(\pi) = \mathrm{crank}(\pi_1) + \frac{1}{2}\mathrm{srank}(\pi)
+ \Psi(\pi),
\label{eq:stcrankdef} 
\eeq
where $\pi_1$ is determined by $\pi \bijone (\pi_1,\pi_2)$,
and the correction term $\Psi(\pi)=1$ if $\pi$ is of type B and zero, otherwise.
We note that
\beq
\mathrm{stcrank}(\pi_A) = -1 + \frac{1}{2}\mathrm{srank}(\pi_A),
\label{eq:stcrankprop1} 
\eeq
and
\beq
\mathrm{stcrank}(\pi_B) = 1 + \frac{1}{2}\mathrm{srank}(\pi_B).
\label{eq:stcrankprop2} 
\eeq
Equipped with the definitions above, we can now prove the following lemma.

\begin{lemma}
\label{lemma1}
If
\beqs
g(x,y,q) := \sum_{\pi} q^{\abs{\pi}} x^{\mathrm{stcrank}(\pi)}
             y^{\mathrm{srank}(\pi)},
\eeqs
then $g(x,y,q)$ has the product representation
\beqs
g(x,y,q) = \frac{ (q^4;q^4)_\infty (-q;q^2)_\infty }
                {(q^4x,q^4/x, q^2y^2x, q^2/(y^2x); q^4)_\infty},
\eeqs
where
\beqs
(a_1,a_2,a_2,\dots ;q)_\infty = (a_1;q)_\infty (a_2;q)_\infty (a_3;q)_\infty
\cdots.
\eeqs
\end{lemma}
\begin{proof}
If $\pi$ is not of type B and $\pi \bijone (\pi_1,\pi_2)$, then
using (\ref{eq:bij1prop1})-(\ref{eq:bij1prop2}), (\ref{eq:stcrankdef})
we find that
\beq
q^{\abs{\pi}} x^{\mathrm{stcrank}(\pi)} y^{\mathrm{srank}(\pi)}
= q^{4\abs{\pi_1} + \abs{\pi_2}} x^{\mathrm{crank}(\pi_1)}
  (xy^2)^{\mathrm{srank}(\pi_2)/2}.
\label{eq:qxy} 
\eeq
On the other hand, if $\pi=\pi_A$ and $\pi_A \bijtwo \pi_B$, then
\begin{align}
&q^{\abs{\pi_A}} x^{\mathrm{stcrank}(\pi_A)} y^{\mathrm{srank}(\pi_A)}
+
q^{\abs{\pi_B}} x^{\mathrm{stcrank}(\pi_B)} y^{\mathrm{srank}(\pi_B)}
\nonumber\\
&\quad=
q^{\abs{\pi_A}} (x+x^{-1}-1) (xy^2)^{\mathrm{srank}(\pi_A)/2}
+
q^{\abs{\pi_B}} x^0 (xy^2)^{\mathrm{srank}(\pi_B)/2}.
\label{eq:qxy2} 
\end{align}
Here we have used (\ref{eq:bij1prop1})-(\ref{eq:bij1prop2})
and (\ref{eq:bij2prop1})-(\ref{eq:stcrankprop2}).

Equations (\ref{eq:qxy}), (\ref{eq:qxy2}) imply that
\beq
\sum_{\pi} q^{\abs{\pi}} x^{\mathrm{stcrank}(\pi)} y^{\mathrm{srank}(\pi)}
=
\sum_{\pi_1} q^{4\abs{\pi_1}} w(x,\pi_1)
\sum_{\pi_2} q^{\abs{\pi_2}} (xy^2)^{\mathrm{srank}(\pi_2)/2}
\label{eq:stcrankgfid} 
\eeq
where
\beq
w(x,\pi_1) =
\begin{cases}
x+x^{-1}-1, &\mbox{if $\pi_1=(1)$,}\\
x^{\mathrm{crank}(\pi_1)}, &\mbox{otherwise.}
\end{cases}
\label{eq:wdef} 
\eeq
We note that in the first sum on the right side of (\ref{eq:stcrankgfid})
the summation is over unrestricted partitions $\pi_1$, and in the
second sum the summation is over partitions $\pi_2$ with no repeated even parts.
Finally, recalling (\ref{eq:crankgf}) with $q\to q^4$ and (\ref{eq:srankprodid})
with $y^2\to xy^2$, we obtain
\beqs
\sum_{\pi} q^{\abs{\pi}} x^{\mathrm{stcrank}(\pi)} y^{\mathrm{srank}(\pi)}
=
\frac{ (q^4;q^4)_\infty }{ (xq^4,q^4/x;q^4)_\infty }
\cdot
\frac{ (-q;q^2)_\infty }{ (xy^2q^2,q^2/(xy^2);q^4)_\infty },
\eeqs
as desired.
\end{proof}

Next we show that
\begin{align}
&\mbox{the coefficient of $q^{5n+4}$ in $g(\xi,1,q)=0$},
\label{eq:coeffz1}\\ 
&\mbox{the coefficient of $q^{5n+4}$ in $g(\xi,\sqrt{-1},q)=0$},
\label{eq:coeffzi} 
\end{align}
where $\xi$ is a primitive fifth root of unity ($\xi^5=1$).
We use the method of \cite{G1}.
We need Jabobi's triple product identity 
\beq
\sum_{n=-\infty}^\infty z^n q^{n^2} = (q^2,-qz-q/z;q^2)_\infty,
\label{eq:jtp} 
\eeq
which implies that
\beq
(q^4;q^4)_\infty (-q;q^2)_\infty = (q^4,-q^3,-q;q^4)_\infty
= \sum_{n=-\infty}^\infty q^{2n^2+n} = \sum_{k\ge0}q^{T_k},
\label{eq:jtpa} 
\eeq
and
\beq
(q^2\xi^2,q^2/\xi^2,q^2;q^2)_\infty = 
\frac{1}{1-\xi^2} \sum_{m\ge0} (-1)^m q^{2T_m} \xi^{-2m}(1-\xi^{4m+2}).
\label{eq:jtpb} 
\eeq
Here $T_k = k(k+1)/2$. 
By Lemma \ref{lemma1} and equations
(\ref{eq:jtpa}) and (\ref{eq:jtpb}) we have
\begin{align}
&g(\xi,1,q) = \frac{ \sum_{k\ge0}q^{T_k} }
                   {(q^4\xi,q^4/\xi,q^2\xi,q^2/\xi;q^4)_\infty}
=
\frac{(q^2\xi^2,q^2/\xi^2,q^2;q^2)_\infty}
     { (q^{10};q^{10})_\infty}
\sum_{k\ge0}q^{T_k} \nonumber \\ 
&\quad
= \frac{1}{1-\xi^2} \frac{1}{(q^{10};q^{10})_\infty}
\sum_{k,m\ge0} (-1)^m q^{2T_m+T_k} \xi^{-2m}(1-\xi^{4m+2}).
\label{eq:gz1id} 
\end{align}
Note that $2T_m + T_k \equiv 4 \pmod{5}$ iff $k\equiv m\equiv 2\pmod{5}$,
but then $1-\xi^{4m+2}=0$. This proves (\ref{eq:coeffz1}).
The proof of (\ref{eq:coeffzi}) is analogous.

Let $P_i(k,m,n)$ denote the number of partitions of $n$ with
$\mathrm{srank}\equiv i\pmod{4}$ and $\mathrm{stcrank}\equiv k\pmod{m}$.
Clearly,
\begin{align}
\sum_{k=0}^4 \xi^k \sum_{n\ge0}P_0(k,5,n) q^n &=
\frac{ g(\xi,1,q) + g(\xi,\sqrt{-1},q)}{2},\label{eq:P0gf} \\ 
\sum_{k=0}^4 \xi^k \sum_{n\ge0}P_2(k,5,n) q^n &=
\frac{ g(\xi,1,q) - g(\xi,\sqrt{-1},q)}{2},\label{eq:P2gf}  
\end{align}
Combining (\ref{eq:coeffz1})-(\ref{eq:coeffzi}) and
(\ref{eq:P0gf})-(\ref{eq:P2gf}) we find that
\beq
\sum_{k=0}^4 \xi^k P_i(k,5,5n+4)=0, (\mbox{for $i=0,2$}),
\label{eq:Piz} 
\eeq
which implies that
\beq
P_i(0,5,5n+4) = P_i(1,5,5n+4) = \cdots = P_i(4,5,5n+4).
\label{eq:Pi04} 
\eeq
On the other hand
\beq
p_i(5n+4)=\sum_{k=0}^4 P_i(k,5,5n+4),
\label{eq:pi54} 
\eeq
so that
\beq
P_i(k,5,5n+4) = \frac{1}{5} p_i(5n+4),
\label{eq:mainresult} 
\eeq
for $i=0,2$ and $k=0,1,2,3,4$.
Thus, we have proved the main result of this section.
\begin{theorem}
\label{theorem1}
The residue of the partition statistic stcrank mod $5$ divides the
partitions enumerated by $p_i(5n+4)$ with $i=0,2$ into
five equinumerous classes.
\end{theorem}
We illustrate this theorem in Table 1 below for the $30$
partitions of $9$. These partitions are organized into five classes
with six members each. In each class the first 4 members have
$\textrm{srank}\equiv0\pmod{4}$ and the remaining two members have
$\textrm{srank}\equiv2\pmod{4}$.

%

\begin{table}[ht]
\caption{}\label{table1}
\beqs
\renewcommand\arraystretch{1.1}
\begin{array}{|c|c|c|c|c|c|}
\hline
&\mathrm{stcrank}\equiv0\mymod{5} &
1\mymod{5} &
2\mymod{5} &
3\mymod{5} &
4\mymod{5} \\
\hline
\mathrm{srank}\equiv0 &({3}^{3}) &
({1}^{5},{2}^{2}) &
({1}^{4},{2}^{1},{3}^{1}) &
({1}^{1},{2}^{4}) &
({1}^{9}) \\
\mymod{4} &({1}^{3},{2}^{1},{4}^{1}) &
({1}^{4},{5}^{1}) &
({1}^{3},{3}^{2}) &
({1}^{6},{3}^{1}) &
({1}^{2},{2}^{2},{3}^{1}) \\
& ({1}^{1},{3}^{1},{5}^{1}) &
({1}^{2},{2}^{1},{5}^{1}) &
({1}^{1},{4}^{2}) &
({1}^{1},{2}^{1},{6}^{1}) &
({2}^{3},{3}^{1}) \\
 & ({4}^{1},{5}^{1}) &
({9}^{1}) &
({2}^{2},{5}^{1}) &
({2}^{1},{7}^{1}) &
({1}^{2},{7}^{1}) \\                
\hline
\mathrm{srank}\equiv2 & ({1}^{3},{2}^{3}) &
({1}^{1},{2}^{1},{3}^{2}) &
({1}^{5},{4}^{1}) &
({1}^{7},{2}^{1}) &
({2}^{1},{3}^{1},{4}^{1}) \\
\mymod{4} &({1}^{3},{6}^{1}) &
({1}^{2},{3}^{1},{4}^{1}) &
({1}^{1},{8}^{1}) &
({1}^{1},{2}^{2},{4}^{1}) &
({3}^{1},{6}^{1}) \\
\hline
\end{array}
\eeqs
\end{table}

Finally, we note that the equation
\beq
\mathrm{srank}(\pi) = - \mathrm{srank}(\pi') \label{eq:srconjprop} 
\eeq
implies that a partition $\pi$ is self-conjugate only if $\mathrm{srank}(\pi)=0$.
This means that the involution $\pi \longrightarrow \pi'$ has no fixed points
if $\mathrm{srank}(\pi)\equiv2\pmod{4}$. Hence, $2\mid p_2(5n+4)$ and
by (\ref{eq:andrefine}) we have the stronger congruence
\beqs
p_2(5n+4) \equiv 0 \pmod{10}.
\eeqs

\section{$t$-cores}
In this section we recall some basic facts about $t$-cores and briefly
review the development in \cite{GKS}. A partition $\pi$ is a called
a $t$-core, if it has no rim hooks of length $t$ \cite{JK}.
We let $a_t(n)$ denote the number of partitions of $n$ which are $t$-cores.
In what follows, $\pi_{\mbox{$t$-core}}$ denotes a $t$-core partition.
Given the diagram of a partition $\pi$ we label a cell
in the $i$-th row and $j$-th column by the least nonnegative integer congruent
to $j-i\pmod{t}$. The resulting diagram
is called a $t$-residue diagram \cite[p.84]{JK}. 

Let $P$ be the set of all partitions and $P_{\mbox{$t$-core}}$ be the
set of all $t$-cores. There is well-known bijection which goes back to
Littlewood \cite{L}.
$\phi_1 \,:\, P \rightarrow P_{\mbox{$t$-core}} \times P \times \cdots \times P$,
\begin{align}
\phi_1(\pi) &= (\pi_{\mbox{$t$-core}}, \vec{\pihat_t}), \label{eq:phi1} \\ 
\vec{\pihat_t} &= (\pihat_0, \pihat_1, \pihat_2, \dots, \pihat_{t-1}),
\label{eq:pihatvecdef} 
\end{align}
such that
\beq
\abs{\pi} = \abs{\pi_{\mbox{$t$-core}}} 
+ t \sum_{i=0}^{t-1}\abs{\pihat_i}.
\label{eq:phi1a} 
\eeq
This bijection is described in more detail in \cite{JK}, \cite{GKS} and \cite{G3}.
The following identity is an immediate corollary of this bijection.
\beq
\frac{1}{(q)_\infty} = \sum_{n\ge0}p(n) q^n
= \frac{1}{(q^t;q^t)_\infty^t} \sum_{n\ge0} a_t(n) q^n.
\label{eq:phi1cor} 
\eeq
It can be rewritten as
\beq
\sum_{n\ge0} a_t(n) q^n = \frac{(q^t;q^t)_\infty^t}{(q)_\infty}.
\label{eq:tcoregfid} 
\eeq

There is another bijection $\phi_2$, introduced in \cite{GKS}.
It is for $t$-cores only.
$\phi_2\,:\, P_{\mbox{$t$-core}} \rightarrow \{\vec{n}
=(n_0, n_1, \dots, n_{t-1}) \,:\, n_i\in\Z, n_0+\cdots+n_{t-1}=0\}$, 
\beq
\phi_2(\pi_{\mbox{$t$-core}}) = \vec{n}=(n_0,n_1,n_2, \dots, n_{t-1}).
\label{eq:phi2} 
\eeq
We call $\vec{n}$ an $n$-vector. It has the following properties.
\beq
\vec{n}\in\Z^t,\qquad \vec{n}\cdot \vec{1}_t = 0,
\label{eq:phi2prop1} 
\eeq
and
\beq
\abs{\pi_{\mbox{$t$-core}}}= \frac{t}{2} \sum_{i=0}^{t-1} n_i^2
+ \sum_{u=0}^{t-1} i n_i,
\label{eq:phi2prop2} 
\eeq
where the $t$-dimensional vector $\vec{1}_t$ has all components
equal to $1$. The generating function identity that corresponds to
this second bijection is

\beq
\sum_{n\ge0} a_t(n) q^n = 
\sum_{\substack{\vec{n}\in\Z^t \\ \vec{n}\cdot\vec{1}_t=0}}
q^{\frac{t}{2}||\vec{n}||^2 + \vec{b}_t\cdot\vec{n}}.
\label{eq:tcoregfid2} 
\eeq

Here
\beq
||\vec{n}||^2=\sum_{i=0}^{t-1} n_i^2,\quad\mbox{and}\quad
\vec{b}_t=(0,1,2, \dots, t-1).
\label{eq:nvecdef} 
\eeq
To construct the $n$-vector of $\pi_{\mbox{$t$-core}}$ in 
(\ref{eq:phi2}), we follow \cite{G3} and define
\beq
\vec{r}(\pi_{\mbox{$t$-core}}) = (r_0, r_1, r_2, \dots, r_{t-1}),
\label{eq:rvec} 
\eeq
where for $0\le i\le t-1$, $r_i(\pi_{\mbox{$t$-core}})$ denotes the number
of cells labelled $i\pmod{t}$ in the $t$-residue diagram of 
$\pi_{\mbox{$t$-core}}$.  Then (\ref{eq:phi2}) can be given explicitly
as
\beq
\phi_2(\pi_{\mbox{$t$-core}}) = \vec{n} = 
(r_0-r_1, r_1-r_2, r_2-r_3, \dots, r_{t-1}-r_0).
\label{eq:rvec2} 
\eeq
We note that $\tfrac{t}{2}||\vec{n}||^2$ is a multiple of $t$ since
$\vec{n}\cdot\vec{1}_t=0$. Hence by (\ref{eq:phi1cor}) and (\ref{eq:tcoregfid2}) 
we have
\beq
\sum_{n\ge0} a_t(tn+\delta) q^{tn+\delta}
=
\sum_{\substack{\vec{n}\in\Z^t,\ \vec{n}\cdot\vec{1}_t=0 \\
\vec{n}\cdot\vec{b}_t\equiv \delta\pmod{t}}}
q^{\frac{t}{2}||\vec{n}||^2 + \vec{b}_t\cdot\vec{n}},
\label{eq:tcoresift} 
\eeq
and
\beq
\sum_{n\ge0} p(tn+\delta) q^{n}
=
\frac{1}{(q)_\infty^t}
\sum_{n\ge0} a_t(tn+\delta) q^{n},
\label{eq:psift} 
\eeq
where
$\delta=0,1,2$, \dots, $t-1$.

We now assume $t=5$. For the case $\delta=4$ the right side of 
(\ref{eq:tcoresift}) can be simplified using the the following change
of variables.
\begin{align}
n_0 &=\alpha_0 + \alpha_4, \nonumber \\
n_1 &=-\alpha_0 + \alpha_1 + \alpha_4, \nonumber \\
n_2 &=-\alpha_1 + \alpha_2, \label{eq:ntoa} \\ 
n_3 &=-\alpha_2 + \alpha_3 - \alpha_4, \nonumber \\
n_4 &=-\alpha_3 - \alpha_4, \nonumber 
\end{align}
We find $\vec{n}$ is an $n$-vector satisfying 
$\vec{n}\cdot\vec{b}_5\equiv4\pmod{5}$ if and only if
\beq
\vec{\alpha} = (\alpha_0, \alpha_1, \alpha_2, \alpha_3, \alpha_4)\in\Z^5
\label{eq:avec1} 
\eeq
and
\beq
\alpha_0 + \alpha_1 + \alpha_2 + \alpha_3 + \alpha_4 = 1.
\label{eq:avec2} 
\eeq
We call $\vec{\alpha}$ and $\alpha$-vector.
Hence, by (\ref{eq:tcoresift}) and (\ref{eq:psift}) we have
\beq
\sum_{n\ge0} a_5(5n+4) q^{n+1}
=
\sum_{\substack{\vec{\alpha}\cdot\vec{1}_5=1 \\
\vec{\alpha}\in\Z^5}}
q^{Q(\vec{\alpha})},
\label{eq:5coresift} 
\eeq
and
\beq
\sum_{n\ge0} p(5n+4) q^{n+1}
=
\frac{1}{(q)_\infty^5}
\sum_{\substack{\vec{\alpha}\cdot\vec{1}_5=1 \\
\vec{\alpha}\in\Z^5}}
q^{Q(\vec{\alpha})},
\label{eq:psift5} 
\eeq
where
\beq
Q(\vec{\alpha}) = ||\vec{\alpha}||^2 - 
                (\alpha_0\alpha_1 + \alpha_1\alpha_2 + \cdots + \alpha_4\alpha_0),
\label{eq:Qadef} 
\eeq
If $\abs{\pi}\equiv4\pmod{5}$ and $t=5$, we can combine 
bijections $\phi_1$ and $\phi_2$ into a single bijection
\beq
\Phi(\pi) = (\vec{\alpha},\vec{\pihat_5}),
\label{eq:Phi} 
\eeq
such that
\beq
\abs{\pi} = 5Q(\vec{\alpha}) - 1 + 5 \sum_{i=0}^4 \abs{\pihat_i}.
\label{eq:Phi1} 
\eeq
Next, following \cite{GKS} we define the $5$-core crank of $\pi$ when
$\abs{\pi}\equiv4\pmod{5}$ as
\beq
c_5(\pi)
=1 + \sum_{i=0}^4 i \alpha_i
\equiv 2(1 + n_0 - n_1 - n_2 + n_3)
\equiv 2 + \sum_{i=-2}^2 i r_{2-i} \pmod{5},
\label{eq:5ccrank} 
\eeq
where $\alpha$ is determined by (\ref{eq:Phi}).

It is easy to check that $Q(\vec{\alpha})$ in ({\ref{eq:Phi1}) remains invariant 
under the following cyclic permutation
\beq
\widehat{C}_1(\vec{\alpha}) = (\alpha_4,\alpha_0,\alpha_1,\alpha_2,\alpha_3),
\label{eq:cycperm} 
\eeq
while $c_5(\pi)$ increases by $1\pmod{5}$ under the map
\beq
\widehat{O}(\pi) = \Phi^{-1}(\widehat{C}_1(\vec{\alpha}), \vec{\pihat}_5).
\label{eq:perm} 
\eeq
In other words, if $\abs{\pi}\equiv4\pmod{5}$, then
\beq
\abs{\pi} = \abs{\widehat{O}(\pi)},
\label{eq:permprop1} 
\eeq
and
\beq
c_5(\pi)+1\equiv c_5(\widehat{O}(\pi))\pmod{5}.
\label{eq:permprop2} 
\eeq
This suggests that all partitions of $5n+4$ can be organized into orbits.
Each orbit consists of five distinct members:
\beq
\pi,\ \widehat{O}(\pi),\ \widehat{O}^2(\pi),\ \widehat{O}^3(\pi),\ 
\widehat{O}^4(\pi),
\label{eq:orb} 
\eeq
and each element of the orbit has a distinct $5$-core crank (mod $5$).
Clearly, the total number of such orbits is $\tfrac{1}{5} p(5n+4)$,
and so $p(5n+4)\equiv0\pmod{5}$. This summarizes the combinatorial
proof of ({\ref{eq:ram5}) given in \cite{GKS}.
If we apply the map $\widehat{O}$ (\ref{eq:perm}) to the partitions of 
$5n+4$ that are $5$-cores, we find that
\beq
a_5^0(5n+4)= a_5^1(5n+4)= \cdots = a_5^4(5n+4),
\label{eq:5corerels} 
\eeq
where, for $0\le j\le4$,  $a_5^j(n)$ denotes the number of partitions of $n$
that are $5$-cores with $5$-core crank congruent to $j$ modulo $5$.
Hence,
\beq
a_5^j(5n+4)= \frac{1}{5} a_5(5n+4),\quad j=0,1,\dots,4,
\label{eq:5corerels2} 
\eeq
which proves that
\beq
a_5(5n+4)\equiv0\pmod{5}.
\label{eq:5corecong} 
\eeq
Actually, more is true. We have 
\beq
a_5(5n+4) = 5a_5(n).
\label{eq:5corerel} 
\eeq
We sketch the combinatorial proof of (\ref{eq:5corerel}) given in \cite{GKS}.
See also \cite{G3}. The map
$\theta:\,P_{\mbox{$5$-core}}(n) \longrightarrow P_{\mbox{$5$-core}}^0(5n+4)$,
defined in terms of $n$-vectors as
\begin{align}
\vec{n} \mapsto \vec{n}'
&=(n_1+2n_2+2n_4+1,
-n_1-n_2+n_3+n_4+1,
2n_1+n_2+2n_3, \nonumber \\
&\qquad -2n_2-2n_3-n_4-1,
-2n_1-n_3-2n_4-1),
\label{eq:theta} 
\end{align}
is a bijection. Here $P_{\mbox{$5$-core}}(n)$ is the set of all $5$-cores of $n$,
and $P_{\mbox{$5$-core}}^0(n)$ is set of all $5$-cores of $n$ with $5$-core
crank congruent to zero modulo $5$. Since $\theta$ is a bijection, we have
\beq
a_5(n)=a_5^0(5n+4).
\label{eq:5corerel2} 
\eeq
The proof of (\ref{eq:5corerel}) easily follows from 
(\ref{eq:5corerels2}) and (\ref{eq:5corerel2}).
Finally, we remark that Ramanujan's result (\ref{eq:rambest}) is
a straightforward consequence of (\ref{eq:psift}) with$(t,\delta)=(5,4)$,
(\ref{eq:5corerel}), and (\ref{eq:tcoregfid}) with $t=5$.

\section{Refinement of Ramanujan's mod $5$ congruence, the srank and
the $5$-core crank}

In the previous section we discussed the combinatorial proof in \cite{GKS} of
Ramanujan's congruence (\ref{eq:ram5}) using the the $5$-core crank 
(\ref{eq:5ccrank}). It is somewhat unexpected that the $5$-core crank can
be employed to prove the refinement (\ref{eq:andrefine}) as well.

In fact, we were amazed to discover the following elegant formulas
\begin{align}
\mathrm{srank}(\pi_{\mbox{$5$-core}}) &\equiv \sum_{i=0}^4 (n_i + i)^3 \pmod{4}, 
\label{eq:elegant1} \\ 
\mathrm{srank}(\pi) &\equiv \mathrm{srank}(\pi_{\mbox{$5$-core}}) 
    + \sum_{i=0}^4 \mathrm{srank}(\pihat_i) \label{eq:elegant2} \\ 
&\quad + 2 \sum_{i=0}^4 \abs{\pihat_i}(n_i + i) \pmod{4},
\nonumber
\end{align}
where $\pi_{\mbox{$5$-core}}$, 
$\vec{\pihat}=(\pihat_0, \pihat_1, \pihat_2, \pihat_3, \pihat_4)$ are determined 
by (\ref{eq:phi1})
with $t=5$, and
\beqs
\vec{n}=(n_0,n_1,\dots,n_4)=\phi_2(\pi_{\mbox{$5$-core}}).
\eeqs
In spite of their simple appearance, the above formulas are far from obvious.
The proof of (\ref{eq:elegant1})-(\ref{eq:elegant2}) will be given elsewhere.
Here we restrict our attention to some implications of 
(\ref{eq:elegant1})-(\ref{eq:elegant2}).

First, we note that if $\abs{\pi_{\mbox{$5$-core}}}\equiv4\pmod{5}$, then
(\ref{eq:elegant1}) can be written in terms an $\alpha$-vector (\ref{eq:ntoa}) as
\beq
\mathrm{srank}(\pi_{\mbox{$5$-core}})
\equiv \alpha_0\alpha_1(\alpha_0-\alpha_1) + \alpha_1\alpha_2(\alpha_1-\alpha_2)
+ \cdots +
\alpha_4\alpha_0(\alpha_4-\alpha_0) \pmod{4}.
\label{eq:sravec} 
\eeq
Similarly, if $\abs{\pi}\equiv4\pmod{5}$, then
\begin{align}
\mathrm{srank}(\pi) &\equiv
\alpha_0\alpha_1(\alpha_0-\alpha_1) + \cdots +
\alpha_4\alpha_0(\alpha_4-\alpha_0)
+ \sum_{i=0}^4 \mathrm{srank}(\pihat_i) \nonumber \\
& \quad + 2\{ (\alpha_0+\alpha_4)\abs{\pihat_0} + (\alpha_2+\alpha_3)\abs{\pihat_1}
+ (\alpha_1+\alpha_2)\abs{\pihat_2} \label{eq:sravec2} \\ 
& \quad + (\alpha_0+\alpha_1)\abs{\pihat_3} + (\alpha_3+\alpha_4)\abs{\pihat_4}\}
\pmod{4}.
\nonumber
\end{align}
Remarkably, (\ref{eq:sravec}) suggests that $\mathrm{srank}(\pi_{\mbox{$5$-core}})$
with $\abs{\pi_{\mbox{$5$-core}}}\equiv4\pmod{5}$ remains invariant
mod $4$ under the cyclic permutation (\ref{eq:cycperm}), and we have
the following refinement of (\ref{eq:5corerels2}):
\beq
a_{5,i}^j(5n+4)= \frac{1}{5} a_{5,i}(5n+4), 
\label{eq:5corerelrefine} 
\eeq
where $j=0$,\dots,$4$ and $i=0$, $2$. Here $a_{5,i}(n)$ denotes the number of 
$5$-cores of $n$ with $\mathrm{srank}\equiv i\pmod{4}$, and
$a_{5,i}^j(n)$ denotes the number of
$5$-cores of $n$ with $\mathrm{srank}\equiv i\pmod{4}$ and $5$-core crank 
$\equiv j\pmod{5}$. Moreover, it is not difficult to verify that the map $\theta$,
given by (\ref{eq:theta}), preserves the srank mod $4$. Indeed, recalling
that $n_0+n_1+n_2+n_3+n_4=0$ we find after some simplication that
\begin{align}
\sum_{i=0}^4 ( (n_i+i)^3 - (n_i'+i)^3) 
&\equiv 2( n_0n_2(n_0+n_2) + n_1n_3(n_1+n_3)+ n_2n_3(n_2+n_3) \nonumber \\
& \quad + n_1(n_1+1) + n_2(n_2+1) + n_3(n_3+1)) \label{eq:invarmod4} \\ 
& \equiv 0 \pmod{4},
\nonumber
\end{align}
where $\vec{n}'$ is defined in (\ref{eq:theta}). Hence, (\ref{eq:5corerel2}) and 
(\ref{eq:5corerel})
can be refined as
\beq
a_{5,i}(n)=a_{5,i}^0(5n+4), \qquad (i=0,2), \label{eq:5corerelrefine2} 
\eeq
and
\beq
a_{5,i}(5n+4)=5 a_{5,i}(n), \qquad (i=0,2), \label{eq:5corerelrefine3} 
\eeq
respectively.

It is less trivial to prove the $5$-core crank analogue of 
Theorem \ref{theorem1}. Namely,
\begin{theorem}
\label{theorem2}
The residue of the $5$-core crank  mod $5$ divides the
partitions enumerated by $p_i(5n+4)$ with $i=0,2$ into
five equal classes.
\end{theorem}
\begin{proof}
We sketch a proof assuming (\ref{eq:sravec2}) and (\ref{eq:5corerelrefine}) hold.
We define the cyclic shift operator $\widehat{C}_2$ by
\beq
\widehat{C}_2(\vec{\pihat}_5) = (\pihat_4,\pihat_2,\pihat_3,\pihat_0,\pihat_1),
\label{eq:C2def} 
\eeq
Next, we use (\ref{eq:C2def}) to modify (\ref{eq:perm}) as
\beq
\widehat{O}_s(\pi) = \Phi^{-1}(\widehat{C}_1(\vec{\alpha}), 
                               \widehat{C}_2(\vec{\pihat}_5)),
\label{eq:newperm} 
\eeq
where
$\Phi(\pi)=(\vec{\alpha},\vec{\pihat}_5)$. Fix $i=0,2$.
By (\ref{eq:sravec2}) we see that $\widehat{O}_s$ preserves the srank mod $4$,
and we may assemble all partitions of $5n+4$ with
$\mathrm{srank}\equiv i\pmod{4}$ into orbits:
\beqs
\pi,\ \widehat{O}_s(\pi),\ \widehat{O}_s^2(\pi),\ \widehat{O}_s^3(\pi),\
\widehat{O}_s^4(\pi),
\eeqs
where $\pi$ is some partition of $5n+4$ with $\mathrm{srank}(\pi)\equiv i\pmod{4}$.
As before, each orbit contains exactly five members and the $5$-core
crank increases by $1$ mod $5$ along the orbit. The number of these
orbits is $\tfrac{1}{5} p_i(5n+4)$, consequently $p_i(5n+4)\equiv0\pmod{5}$ and 
the result follows.
\end{proof}

Theorem \ref{theorem2} is illustrated below in Table \ref{table2},
which contains all $30$ partitions of $9$, organized into $6$ orbits.
Each row in this table represents an orbit, and the first row lists
all partitions of $9$ that are $5$-cores.

\begin{table}[ht]
\caption{}\label{table2}
\beqs
\renewcommand\arraystretch{1.1}
\begin{array}{|c|c|c|c|c|c|}
\hline
&\mbox{$5$-core crank}\equiv0\mymod{5} &
1\mymod{5} &
2\mymod{5} &
3\mymod{5} &
4\mymod{5} \\
\hline
\mathrm{srank}\equiv0 & ({1}^{4},{5}^{1}) & 
({1}^{3},{3}^{2}) & 
({1}^{4},{2}^{1},{3}^{1}) & 
({1}^{1},{2}^{1},{6}^{1}) & 
({2}^{2},{5}^{1}) \\ 
\mymod{4} & ({1}^{5},{2}^{2}) & 
({2}^{3},{3}^{1}) & 
({1}^{2},{7}^{1}) & 
({4}^{1},{5}^{1}) & 
({1}^{3},{2}^{1},{4}^{1}) \\ 
 &({3}^{3}) & 
({1}^{9}) & 
({1}^{1},{3}^{1},{5}^{1}) & 
({1}^{2},{2}^{2},{3}^{1}) & 
({9}^{1}) \\ 
 & ({2}^{1},{7}^{1}) & 
({1}^{2},{2}^{1},{5}^{1}) & 
({1}^{1},{2}^{4}) & 
({1}^{6},{3}^{1}) & 
({1}^{1},{4}^{2}) \\ 
\hline
 \mathrm{srank}\equiv2 & ({1}^{3},{2}^{3}) & 
({1}^{3},{6}^{1}) & 
({2}^{1},{3}^{1},{4}^{1}) & 
({1}^{1},{8}^{1}) & 
({1}^{2},{3}^{1},{4}^{1}) \\ 
\mymod{4} & ({3}^{1},{6}^{1}) & 
({1}^{1},{2}^{2},{4}^{1}) & 
({1}^{7},{2}^{1}) & 
({1}^{1},{2}^{1},{3}^{2}) & 
({1}^{5},{4}^{1}) \\ 
\hline
\end{array}
\eeqs
\end{table}

New, we state some new formulas for $a_{5,0}(n)$:
\begin{align}
a_{5,0}(4n) &= a_5(4n), \label{eq:a50form1} \\ 
a_{5,0}(4n+1) &= a_5(4n+1), \label{eq:a50form2} \\ 
a_{5,0}(4n+2) &= 0, \label{eq:a50form3} \\ 
a_{5,0}(4n+3) &= a_5(n). \label{eq:a50form4}  
\end{align}
Formulas (\ref{eq:a50form1})-(\ref{eq:a50form3}) follow from (\ref{eq:elegant1}).
Formula (\ref{eq:a50form4}) is a consequence of the following bijective map,
defined in terms of $n$-vectors by
\beq
\vec{n} \mapsto \vec{n}'
=(2n_1, 1+2n_4, 2n_2, -1+2n_0, 2n_3).
\label{eq:a50trans} 
\eeq
The important properties of (\ref{eq:a50trans}) are
\beq
\abs{\phi_2^{-1}(\vec{n}')} = 4 \abs{\phi_2^{-1}(\vec{n})} + 3,
\label{eq:transprop1} 
\eeq
and
\beq
\mathrm{srank}(\phi_2^{-1}(\vec{n}')) \equiv 0 \pmod{4}.
\label{eq:transprop2}
\eeq
The details will be given elsewhere.

\section{Concluding remarks}
While the stcrank development in section 2 followed naturally from the Andrews
product (\ref{eq:rsgf}), the $5$-core crank development in the previous
section arose in an unexpected fashion.  One may wonder, if there exists
an additional new partition statistic, closely related to the Dyson rank,
whose residue mod $5$ splits the partitions enumerated by $p_i(5n+4)$
with $i=0,2$ into $5$ equal classes.
Finally, we would like to pose the

\noindent
\textbf{Problem}. Is there an analogue of the Stanley rank, which gives
a refinement for Ramanujan's partitions congruences mod $7$ and mod $11$?

{\it Acknowledgements}.
We would like to thank George Andrews and Krishna Alladi for their genuine
interest and encouragement.

\bibliographystyle{amsplain}

\end{document}